\documentclass[11pt]{amsart}
\usepackage{amsthm}
\usepackage{amssymb}
\usepackage{amsfonts}
\usepackage{dsfont}

\def\p{{\varphi}}
\def\w{{\omega}}

\def\cI{{\mathcal I}}

\def\B{{\mathbf{B}}}

\def\E{{\mathds E}}
\def\P{{\mathds P}}
\def\N{{\mathds N}}
\def\R{{\mathds R}}
\def\1{{\mathds 1}}

\def\dint{\textup{d}}

\begin{document}

\title{Monotonicity of functionals of random polytopes}
\date{}

\author{Mareen Beermann, Matthias Reitzner}
\address{University of Osnabrueck, Department of Mathematics, Albrechtstr. 28a, 49076 Osnabrueck, Germany}
\email{mareen.beermann[at]uni-osnabrueck.de, matthias.reitzner[at]uni-osnabrueck.de}



\maketitle

\centerline{\it Dedicated to Imre B\'ar\'any on the occasion of his 70th birthday.}

\section{Introduction}
Let $n$ random points $X_{1},...,X_{n}$ be chosen independently and according to a given density function $\phi$ in $\R^{d}$. We call the convex hull $P_{n}=[X_{1},...,X_{n}]$ of these points a random polytope. Various properties of these objects have been studied in the last decades, e.g. the number of $j$-dimensional faces and the intrinscic volumes. Classical papers dealt with the expected values of these functionals, see e.g. B\'ar\'any \cite{Bar2}, B\'ar\'any and Buchta \cite{BB}, B\'ar\'any and Larman \cite{BL}, Reitzner \cite{Re7}. More recently, distributional properties have been investigated intensively, e.g. variance estimates,  central limit theorems and and large deviation inequalities, see e.g. 
B\'ar\'any, Fodor, and Vigh \cite{BaFoVi10},
B\'ar\'any and Reitzner \cite{BRe2, BRe1}, 
Calka, Schreiber and Yukich \cite{CaScYu13},
Calka and Yukich \cite{CaYu, CaYu15, CaYu17} , 
Pardon \cite{Par11, Par12} 
Reitzner \cite{Re6, Re8, Re7}, 
Schreiber and Yukich \cite{SchYu},
and Vu \cite{Vu1}.

For all these questions the expectation is a central object. 
We denote by $ \E V_d(P_n)$ the expected volume of the random polytope and by $\E f_j(P_n)$ the expectation of the number of $j$-dimensional faces. 

In this short note we concentrate on monotonicity questions concerning $\E V_d(P_n)$ and $\E f_j(P_n)$ which have been investigated in the last years. For more information on random polytopes and related questions we refer to the survey articles \cite{Hugsurv} and \cite{Re9}.

\section{Monotonicity of the volume with respect to set inclusion}
Let $K \subset \R^d$ be a convex set with nonempty interior and $\phi(\cdot)= V_d(K)^{-1} \1_K(\cdot)$, thus the points $X_i$ are chosen according to the uniform distribution in $K$ and $P_n$ is a random polytope in the interior of $K$. It seems to be immediate that increasing the convex body $K$ should also increase the random polytope and thus its volume.

More precisely, assume that $K, L$ are two $d$-dimensional convex  sets. Choose independent uniform random points 
$X_1, \dots , X_n$ in $K$
and
$Y_1, \dots, Y_n$ in $L$, 
and denote by $P_n$ the convex hull $[X_1, \dots, X_n]$, and by $Q_n$ the convex hull $[Y_1, \dots, Y_n]$. 
Is it true that $K \subset L$ implies
\begin{equation}\label{eq:mon-d-gen}
 \E V_d(P_n) \leq  \E V_d(Q_n)?  
\end{equation}
The starting point for the investigation should be a check for the first nontrivial case 
$n=d+1$ where a random simplex is chosen in $K$, resp. $L$. In this form, the 
question was first raised by Meckes \cite{MeckesAMS} in the context of high-dimensional convex 
geometry. 

In dimension one the monotonicity is immediate. In 2012 in a groundbreaking paper by Rademacher \cite{Rademacher} proved
that this is also true in dimension two, but that there are counterexamples for dimensions $d \geq 4$ and $n=d+1$. 
Only recently the three-dimensional case could be settled by Kunis, Reichenwallner and Reitzner \cite{KuReRe17} where monotonicity of $\E V_3(P_4)$ also fails. The question of monotonicity of higher moments $\E V_d(P_n)^k$ was investigated  in \cite{ReichenwallnerReitzner}.

It remains an open problem whether there is a number $N$, maybe depending on $K$ or only on the 
dimension $d$, such that monotonicity holds for $n \geq N$. 

\section{Monotonicity of the number of faces with respect to $n$}
Choose $X_{1},...,X_{n}$ according to a given density function $\phi$ in $\R^{d}$. 
A natural guess is that the expected number $\E f_{j}(P_{n})$ of $j$-dimensional faces behaves monotone if the number of generating points increases. The asymptotic results suggest that at least for random points chosen uniformly in a smooth convex set and or a polytope (see \cite{Bar2}, \cite{BB}, \cite{BL}, \cite{Re7}) the expectation $\E f_j(P_n)$ should be increasing in $n$,
$$ \E f_j(P_n) \leq \E f_j(P_{n+1}) \ \ \ \forall n \in \N .$$
On the other hand B\'ar\'any \cite{Bar1} showed that the behaviour for generic convex sets is extremely complicated and thus monotonicity is not obvious.

The first results concerning this issue have been gained by Devillers et al \cite{DGGMR}. They considered convex hulls of uniformly distributed random points in a convex body $K$. It is proven that for planar convex sets the expected number of vertices $\E f_{0}(P_{n})$ (and thus also edges) is increasing in $n$. Furthermore they showed that for $d\geq 3$ the number of facets $\E f_{d-1}(P_{n})$ is increasing for $n$ large enough if $\lim_{n \rightarrow \infty}\frac{\E f_{d-1}(P_{n}}{An^{c}}=1$ for some constants $A$ and $c>0$, e.g. for $K$ being a smooth convex body.
In the PhD thesis of Beermann \cite{BeDiss} the cases of $\phi$ being the Gaussian distribution or the uniform distribution in a ball are settled. We sketch the proof in the Appendix. The method used for these results was extented by Bonnet et al \cite{Bonnetetal} who settled the cases of random points on the sphere, on a halfsphere, random points chosen according to a certain heavy-tailed distribution, and beta-type distributions.

It should be noted that these results carry over to monotonicity results for convex hulls of random points chosen from a Poisson point process with a suitable density.

All these results only deal with the number of facets. Only in the Gaussian case it seems to be possible to extend this monotonicity results to general $j$-dimensional faces. This is the content of a recent preprint by Kabluchko. But even for other \lq most simple cases\rq\  like uniform points in a ball the general question is widely open.

\section{Appendix: Facet numbers of random polytopes}

Let $\Phi$ be a probability measure in $\R^d$ with density $\p$. Choose $n$ random points $X_1, \dots, X_n$ independently according to $\Phi$, and let $P_n$ be the convex hull of these random points. We start by developing a well known formula for $\E f_{d-1}$. Each $(d-1)$-dimensional face of $P_{n}$ is the convex hull of exactly $d$ random points with probability one. Since $X_{1},...,X_{n}$ are chosen independently and identically it holds
\begin{eqnarray*}
\E f_{d-1}(P_{n}) 
& = & \binom nd \P([X_{1}, \dots ,X_{d}] \text{ is a facet }).
\end{eqnarray*}
We denote by $H_{1, \dots, d}$ the affine hull of the $(d-1)$-dimensional simplex $P_{d}=[x_{1},\dots,x_{d}]$ which divides $\R^{d}$ into the two halfspaces $H_{1, \dots, d}^{+}$ and $H_{1, \dots, d}^{-}$. If $P_{d}$ is a facet then the other points $X_{d+1},\ldots , X_{n}$ are either all located in $H_{1, \dots, d}^{+}$, or in $H_{1, \dots, d}^{-}$. This happens with probability $\Phi(H_{1, \dots, d}^+)^{n-d}$, resp. $\Phi(H_{1, \dots, d}^-)^{n-d}$, hence 
\begin{eqnarray}\label{allgemeinEfacet}
\E f_{d-1}(P_{n}) & = & \binom nd  \int\limits_{\R^{d}} \dots \int\limits_{\R^{d}} 
\big( \Phi(H_{1, \dots, d}^+)^{n-d} + \Phi(H_{1, \dots, d}^-)^{n-d} \big) \, \prod\limits_{i=1}^{d} \phi(x_i) \dint x_i.
\end{eqnarray}
We use the classical affine Blaschke-Petkantschin formula (see \cite{SchnWe3}, Theorem 7.2.1.) to conclude 
\begin{eqnarray*}
\E f_{d-1}(P_{n}) 
& = & 
(d-1)! \binom nd  \int \limits_{S^{d-1}} \int\limits_0^{\infty} \big( \Phi(H(p,\w)^+)^{n-d} + \Phi(H(p,\w)^-)^{n-d} \big) \\
 &  & \hspace*{1.5cm}
\times  \int\limits_{H(p,\w)} \dots \int\limits_{H(p,\w)} \Delta_{d-1}(x_{1},...,x_{d}) \prod\limits_{i=1}^{d} \phi(x_i) \dint x_{i} \  \dint p \dint \omega
\\ & = & 
(d-1)! \binom nd  \int \limits_{S^{d-1}} \int\limits_{-\infty}^\infty \Phi(H(p,\w)^+)^{n-d}
\\ &  & \hspace*{1.5cm}
\int\limits_{H(p,\w)} \dots \int\limits_{H(p,\w)} \Delta_{d-1}(x_{1},...,x_{d}) \prod\limits_{i=1}^{d} \phi(x_i) \dint x_{i} \ \dint p \dint \omega
\end{eqnarray*}
where we parametrize the hyperplane by $H(p,\w)= \{ x\colon \langle x,\omega \rangle =p\}$, and the halfspaces by 
$H(p,\w)^- =\{ x\colon \langle x,\omega \rangle  \leq p\}$ and
$H(p,\w)^+ =\{ x\colon \langle x,\omega \rangle  \geq p\}$.
In the inner integral $\Delta_{d-1}(x_{1},...,x_{d})$ is the $(d-1)$-dimensional volume of the convex hull of $x_{1},...,x_{d}$. We fix the direction $\omega$ and want to prove the montonicity in $n$ of 
\begin{eqnarray*}
\cI(n)
& = & 
\binom nd \int\limits_{\R}\Phi(H(p,\w)^+)^{n-d}
\int\limits_{H(p,\w)} \dots \int\limits_{H(p,\w)} \Delta_{d-1}(x_{1},...,x_{d}) \prod\limits_{i=1}^{d} \phi(x_i) \dint x_{i} \ \dint p 
\end{eqnarray*}
For a given direction $\w$ we put 
$ 
\psi(t) = \int_{H(t, \w) }\phi(x) \dint x \ \text{and}\ 
\Psi(p) = \int_{-\infty}^p \psi(t) \dint t 
$
which defines the push forward measure of $\Phi$ under the projection onto the line $\{t\w\colon t\in \R\}$.
Observe that on the support of $\psi$, the mass of halfspaces $\Psi(p)=\Phi(H(p, \w)^-)$ is an increasing function in $p$ and thus there is an inverse function $\Psi^{-1}(s)$, also increasing, with 
$$ 
\frac{d}{ds} \Psi^{-1} (s) = 
\Big(\frac{d}{dp} \Phi (H(p, \w)^-) \vert_{p=\Psi^{-1} (s)}\Big)^{-1} = 
\Big(\psi (p) \vert_{p=\Psi^{-1} (s)}\Big)^{-1} = 
\Big(\psi (\Psi^{-1} (s)))^{-1} 
$$
for $s \in (0,1)$, and thus $\dint p= (\psi (\Psi^{-1} (s)))^{-1} \dint s$.
Substituting by $s=\Psi(p)=\Phi(H(p,\w)^-)$ we end up with 
\begin{eqnarray*}
\cI(n)
& = & \binom nd \int\limits_{0}^{1} (1-s)^{n-d} \psi(\Psi^{-1}(s))^{d-1} \, \E_{H(\Psi^{-1}(s),\w)} \Delta_{d-1} (X_1, \dots X_d) \ \dint s 
\end{eqnarray*}
where $\E_{H(p,\w)} \Delta_{d-1} (X_1, \dots X_d)$ is the volume of a random simplex where the points $X_1, \dots ,X_d$ are chosen independently according to the normalized density $\p$  in $H(p,\w)$.
Thus to prove monotonicity we have to show that $\triangle_n \cI=\cI(n)-\cI(n-1)$ is positive,
\begin{eqnarray*}\label{eq:pos}
\triangle_n \cI
& = & 
\frac 1n \binom nd \int\limits_{0}^{1} (1-s)^{n-d-1} 
(d-ns ) L(s)^{d-1}\ \dint s 
\end{eqnarray*}
with
\begin{eqnarray*}
L(s) & = & 
\psi(\Psi^{-1}(s)) 
\Big(\E_{H(\Psi^{-1}(s),\w)} \Delta_{d-1} (X_1, \dots X_d)  \Big)^{\frac 1{d-1}}.
\end{eqnarray*}
In the next two sections we will show that in both cases we are interested in, the function $L(s)$ is concave. This is sufficient, because then the graph of $L(s)$ starts at the origin, is above the line $l(s)=L( \frac dn) \frac{ns}{d}$ in $(0, \frac dn)$, meets the line for $s=\frac dn$, and is below the line for $s > \frac dn$.
This yields
\begin{eqnarray*}
\triangle_n \cI
& \geq & 
\frac 1n \binom nd \int\limits_{0}^{1} (1-s)^{n-d-1} 
(d-ns ) l(s)^{d-1}\ \dint s 
\\ & = & 
\frac{n^{d-2}}{d^{d-1}}\binom nd L\Big( \frac dn\Big)^{d-1} 
\underbrace{\int\limits_{0}^{1} (1-s)^{n-d-1}  (d-ns ) s^{d-1}  \ \dint s }_{=d \B(n-d,d)- n \B(n-d,d+1)}
=0\ , 
\end{eqnarray*}
and hence  $\E f_{d-1} (P_n)$ is increasing.

\subsection{The Case of Gaussian Polytopes}
In this case we have $\phi(x)= \frac 1 {(2 \pi)^{d/2}} \exp\{-\sum x_i^2 /2  \}$. 
By the rotation invariance it sufficies to consider the case $\w=(1,0,\dots,0)$ where it is easy to see that $\psi(t)= \frac 1{\sqrt {2\pi}} \exp\{-t^2 /2\}$ and that 
$
\E_{H(p,\w)} \Delta_{d-1} (X_1, \dots X_d) 
$
is independent of $p$. Thus $L(s) = c_d \psi (\Psi^{-1}(s))$.

The continously differentiable function $L(s)$ is concave if and only if its derivative is decreasing. Since
$ \psi'(t) = -t \psi(t)$, it follows that
$$ L'(s) = c_d \frac d{ds} \psi(\Psi^{-1} (s)) = - c_d \Psi^{-1}(s) \psi(\Psi^{-1}(s))  \big(\psi(\Psi^{-1}(s))\big)^{-1}  = - c_d \Psi^{-1} (s).
$$
Clearly, $\Psi(s)$ is increasing in $s$, and therefore $\Phi^{-1} (s)$ too. This implies that $-\Phi^{-1} (s)$ is decreasing and $L(s)$ is concave on $[0,1]$.

\subsection{The Case of Random Polytopes in a Ball}
Assume that $B^d$ is the unit ball of volume $\kappa_d$.
In this case $\phi(x)=\kappa_d^{-1} \1(x \in B^d)$. 
By the rotation invariance it sufficies to consider the case $\w=(1,0,\dots,0)$ where 
$\psi(t)= \kappa_d^{-1} \kappa_{d-1} (1-t^2)^{(d-1)/2} $. 
The intersection of $H(t, \w) \cap B^d$ is always a ball of radius $(1-t^2)^{1/2} $ and the expected volume of a random simplex in $H(t,\w)\cap B^d$ is a constant times  $V_{d-1}(H(t,\w) \cap B^d)=\psi(t)$. The constant is determined explicitly in a paper of Miles \cite{Mi1}. Thus
$
L(s)= c_d \psi(\Psi^{-1}(s))^{\frac d{d-1}} 
$.
We have
$$ L'(s) = 
c_d \frac d{ds} \psi(\Psi^{-1}(s))^{\frac d{d-1}} = 
c_d\, \frac d{d-1}\,  \frac{\psi'(\Psi^{-1}(s))} {\psi(\Psi^{-1}(s))^{\frac {d-2}{d-1}}} =
c_d\, d   \Big(\frac d{dp} \psi(p)^{\frac1{d-1}}\Big) \Big\vert_{p=\Psi^{-1}(s)} .
$$
Because $\psi$ is a concave function, its derivative is decreasing in $p$. Noting that $\Psi^{-1}$ is increasing shows that 
$L'$ is decreasing und thus the function $L(s)$ is concave.


\end{document}